\documentclass[12pt]{article}
\usepackage[cp1251]{inputenc}
\usepackage{amsmath}
\usepackage{amssymb}
\usepackage[english]{babel}
\usepackage[T2A]{fontenc}
\usepackage{epsfig}

\begin{document}

\begin{large}

\centerline{\textbf{Circular tractrices and generalized Dini surfaces}}

\end{large}

\bigskip

\centerline{V.~O.~Gorkavyy\footnote{E-mail: gorkaviy@ilt.kharkov.ua, vasylgorkavyy@gmail.com }  }

\centerline{(B.~Verkin  Institute for Low Temperature Physics and Engineering, Kharkiv, Ukraine)}
\bigskip

{\small {\bf Abstract.} We introduce a particular family of two-dimensional surfaces in
$\mathbb R^4$ which generalize the classical Dini surfaces in $\mathbb R^3$.}

\bigskip

\bigskip

\textbf{1. Introduction}

\medskip

This paper deals with pseudo-spherical surfaces in $n$-dimensional Euclidean spaces $E^n$. It is mainly motivated by the problem settled by Yu. Aminov and A. Sym  whether one can extend the classical theory of Bianchi-Backlund transformations of two-dimensional pseudo-spherical surfaces in $E^3$ to the case of two-dimensional surfaces in $E^n$ with $n\geq 4$ \cite{AmS}.

By definition, a surface in $E^3$ is referred to as {\it pseudo-spherical} if its Gauss curvature is constant negative.
The most known examples of pseudo-spherical surfaces in $E^3$ are named after the famous Italian mathematicians
Eugenio Beltrami and Ulisse Dini. The Beltrami surface (the pseudo-sphere) is a surface of revolution obtained by rotating a particular curve called the tractix. The Dini surfaces are helicoidal surfaces in $\mathbb E^3$ obtained by applying screw rotations to the tractix  \cite{Dini}.

\begin{center}
\epsfig{file=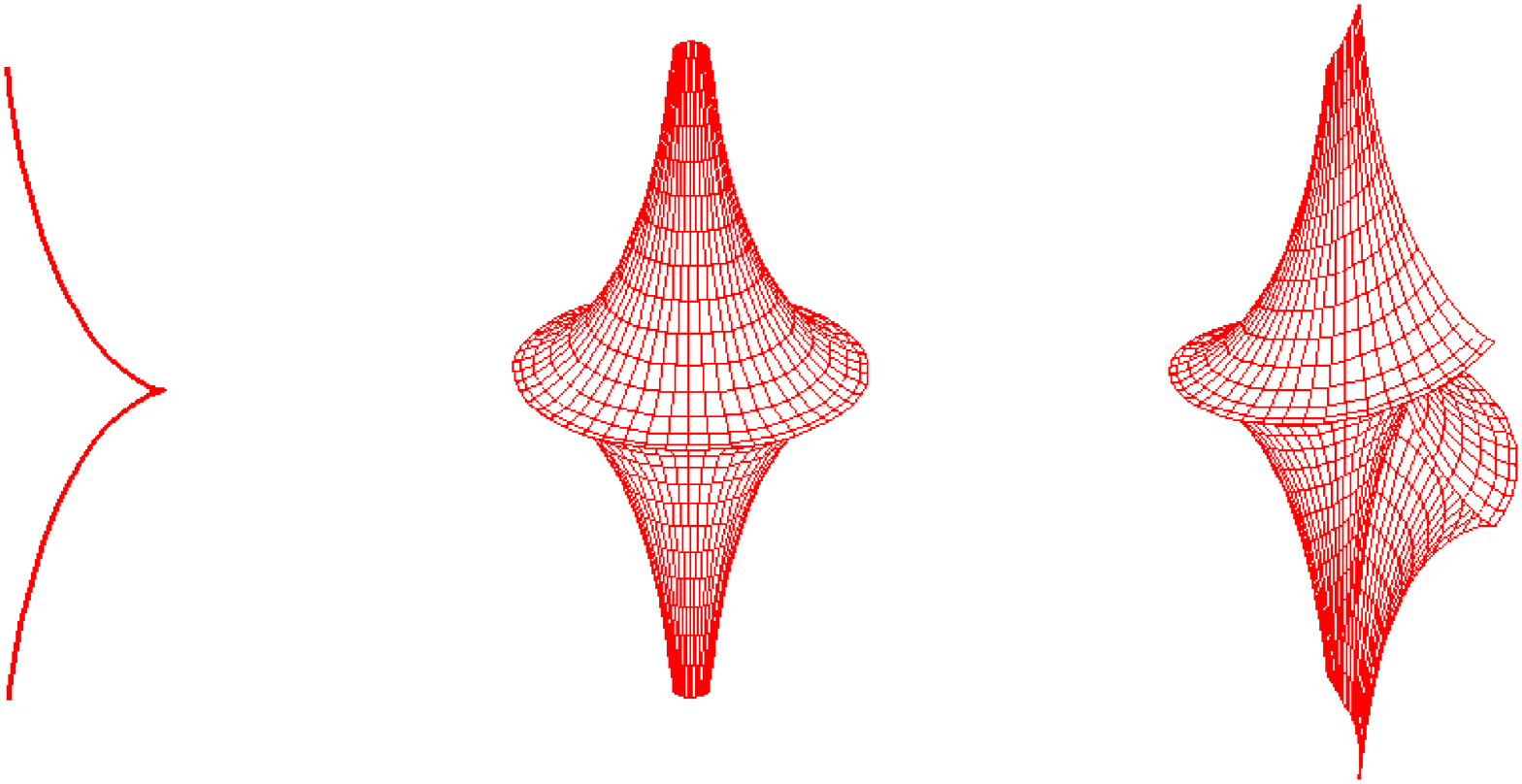, height=40mm, width=120mm}
\end{center}

\smallskip

\centerline{Fig. 1. Tractrix (left), Beltrami surface (center), Dini surface (right)}

\medskip

In frames of the theory of integrable systems pseudo-spherical surfaces in $E^3$ correspond to solutions of the sine-Gordon equation  \cite{Tenenblat}. Particularly, the Beltrami and Dini surfaces correspond to one-soliton solutions of the sine-Gordon equation. Hence they can be characterized as the pseudo-spherical surfaces that admit {\it degenerate} Bianchi-Backlund transformations. Geometrically, this means that if one takes either the Beltrami surface or one of the Dini surfaces and draws appropriate segments of constant length tangent to the meridians (the tractrices) of the surface, then the endpoints of the segments will sweep out not a two-dimensional surface but a straight line, the axe of rotation of the surface.

While exploring the question posed by Yu.Aminov and A. Sym, we was interested in finding surfaces in $E^n$, $n\geq 4$, which inherit geometric properties of the Beltrami and Dini surfaces related to the degeneracy of Bianchi-Backlund transformations. As result, a novel family of pseudo-spherical surfaces in $E^n$, $n\geq 4$, called generalized Beltrami surfaces was described in \cite{G1}, where it was shown that every generalized Beltrami surface admits a degenerate Bianchi transformation like the classical Beltrami surface in $E^3$ do, see also \cite{G2}.

In this note we present another novel family of two-dimensional surfaces in $E^4$ which can pretend to be a geometrically consistent analog of the Dini surfaces in $E^3$.

\newpage

The key ingredient in our studying are {\it circular tractrices} in $E^3$ which we propose to use instead of the classical tractrix. Namely, given a smooth curve $\tilde\gamma$ in $E^3$ and a smooth field of unit\footnote{Instead of unit segments, one can consider segments of an arbitrary fixed length. This more general situation can be reduced to the case we consider in this note by applying dilatations in $\mathbb E^3$.} segments along $\tilde\gamma$, consider a curve $\gamma$ swept out by the endpoints of the unit segments in question. If these segments are tangent to $\gamma$, then $\gamma$ is called a tractrix and $\tilde\gamma$ is called the directrix of $\gamma$. If $\tilde\gamma$ is a circle, then $\gamma$ is called a circular tractrix, c.f. \cite{tr1}, \cite{tr2}.

\begin{center}
\epsfig{file=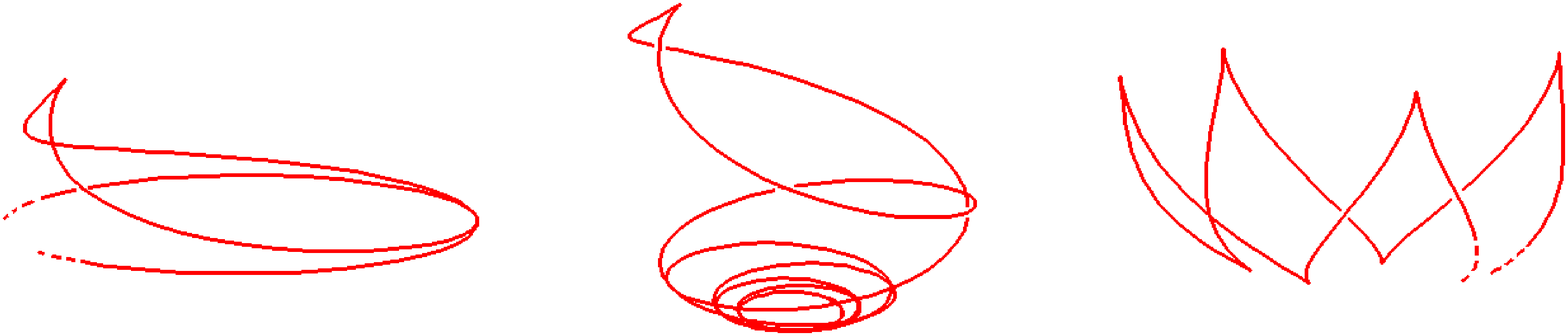, height=30mm, width=90mm}
\end{center}

\smallskip

\centerline{Fig. 2. Examples of circular tractrices in $E^3$}

\medskip

We provide an explicit description for the circular tractrices in $E^3$. Notice that these curves are quite non-trivial and possess a lot of surprising geometrical properties which justify their use as circular analogs of the classical linear tractrix, see \cite{sirosh}. Particularly, for an arbitrary circle $\tilde\gamma$ there exists a one-parametric family of non-congruent circular tractrices whose common directrix is $\tilde\gamma$.

Next, given an arbitrary circular tractrix $\gamma$ in $E^3$, we view $E^3$ as a hyperplane in $E^4$ and split $E^4$ into the direct sum $E^2_1\oplus E^2_2$ so that $E^2_1$ is a two-dimensional subspace of $E^4$ containing the directrix $\tilde\gamma$ of $\gamma$ and $E^2_2$ is its orthogonal complement in $E^4$. Then we apply to  $\gamma$ a one-parametric skew rotation in $E^4$ which is composed of a rotation along the circle $\tilde\gamma$ in $E^2_1$ and a rotation in $E^2_2$. Under this skew rotation, the circular tractrix $\gamma$ sweeps out a two-dimensional helicoidal surface $F^2\subset E^4$.

We analyze fundamental properties of the helicoidal surface $F^2\subset E^4$ and, as the main result, we demonstrate that the Gauss curvature of $F^2$ is constant, $K\equiv K_0$. The value of $K_0$ can be negative, zero or positive, it depends on the radius of the circle $\tilde\gamma$ and on the choice of $\gamma$ in the one-parametric family of circular tractices whose common directrix is $\tilde\gamma$.

Moreover, it is easy to see by construction that if one draws appropriate unit segments tangent to the meridians of the surface $F^2$, then the endpoints of these segments will sweep out the circle $\tilde\gamma$.

Thus, we get a family of helicoidal surfaces of constant Gauss curvature in $E^4$ which inherit geometric properties of the classical Dini surfaces in $E^3$ and hence can be viewed as their circular analogs. We propose to call them {\it the circular Dini surfaces} and hope that these surfaces will be of interest for a further more detailed studying.

\medskip

\bigskip

\textbf{2. General tractrices in $E^n$}

\medskip

Firstly, we describe a general method for constructing a tractrix in $E^n$ with an arbitrarily given directrix, c.f. \cite{tabachnikov}.

Let $\tilde \gamma$ be a smooth curve in $E^n$, $n\geq 2$. Denote by $\tilde f(s)$ the position vector of $\tilde \gamma$ parameterized by an arc length $s$. Take an arbitrary smooth unit vector field along $\tilde\gamma$ represented by a unit vector-function $v(s)$ and consider a curve $\gamma$ in $E^n$ represented by the position vector
\begin{eqnarray}\label{1}
f(s) = \tilde f(s)+v(s).
\end{eqnarray}

By definition, $\gamma$ is a tractrix and $\tilde\gamma$ is its directrix if and only if one has
\begin{eqnarray}\label{2}
[\frac{df}{ds} , v ] \equiv 0,
\end{eqnarray}
where $[,]$ stands for the cross-product of vectors in $E^n$.

\newpage

For reconstructing the tractrix $\gamma$ from its directrix $\tilde\gamma$ we have to substitute (\ref{1}) into (\ref{2}) and then consider (\ref{2}) as a differential equation for the vector-function $v(s)$.

We use the Frenet frame $\xi_1$, ..., $\xi_n$ associated with $\tilde \gamma$ and write down $v$ as a linear combination,
\begin{eqnarray}\label{3}
v = \sum\limits_{j=1}^n v^j \xi_j,
\end{eqnarray}
where $v^j=v^j(s)$, $1\leq j\leq n$  are some functions.

By applying the Frenet-Serre formulae, see \cite{Aminov}, one gets
\begin{eqnarray*}\label{4}
\frac{d v}{ds} = \sum\limits_{j=1}^n \left( \frac{dv^j}{ds} \xi_j + v^j\frac{d\xi_j}{ds} \right)
= (\frac{dv^1}{ds}-v^2k_1) \xi_1 + (\frac{dv^2}{ds}+v^1k_1-v^3k_2 ) \xi_2 + ...\\ ... + (\frac{dv^{n-1}}{ds}+v^{n-2}k_{n-2}-v^nk_{n-1}) \xi_{n-1} + (\frac{dv^n}{ds}+v^{n-1}k_{n-1}) \xi_n ,
\end{eqnarray*}
where $k_1$, $k_2$, ... $k_{n-1}$ are curvatures of $\tilde\gamma\subset E^n$.

Since $\displaystyle{\frac{d\tilde f}{ds} = \xi_1}$, then
\begin{eqnarray}\label{tan}
 \frac{d f}{ds} = \frac{d \tilde f}{ds} + \frac{d v}{ds} = (\frac{dv^1}{ds}-v^2k_1+1) \xi_1 + (\frac{dv^2}{ds}+v^1k_1-v^3k_2 ) \xi_2 + ... \\ ... + (\frac{dv^{n-1}}{ds}+v^{n-2}k_{n-2}-v^nk_{n-1}) \xi_{n-1} + (\frac{dv^n}{ds}+v^{n-1}k_{n-1}) \xi_n.
\end{eqnarray}
Therefore, (\ref{2}) rewrites as follows:
\begin{eqnarray*}
\frac{\displaystyle{\frac{dv^1}{ds}}-v^2k_1+1}{v^1} = \frac{\displaystyle{\frac{dv^2}{ds}}+v^1k_1-v^3k_2 }{v^2} = ...
= \frac{\displaystyle{\frac{dv^{n-1}}{ds}}+v^{n-2}k_{n-2}-v^nk_{n-1}}{v^{n-1}} = \frac{\displaystyle{\frac{dv^n}{ds}}+v^{n-1}k_{n-1}}{v^n} ,
\end{eqnarray*}
or equivalently as follows
\begin{eqnarray}\label{ode_sys}
\frac{d}{ds}
\left(
\begin{array}{l}
v^1 \\ v^2 \\ ... \\ v^{n-1} \\ v^n
\end{array}
\right) =
\left(
\begin{array}{lllllll}
0&k_1& 0 & ...  & 0 & 0 \\
-k_1 & 0 & k_2 & ...  & 0 & 0 \\
...& ... & ... & ...  & ... & ... \\
0 & 0 & 0 & ...  & 0 & k_{n-1} \\
0 & 0 & 0 & ...  &  -k_{n-1} & 0 \\
\end{array}
\right) \cdot
\left(
\begin{array}{l}
v^1 \\ v^2 \\ ... \\ v^{n-1} \\ v^n\\
\end{array}
\right)
-
\left(
\begin{array}{l}
1 \\ 0 \\ ... \\ 0 \\ 0 \\
\end{array} \right) +
\lambda \left( \begin{array}{l}
v^1 \\ v^2 \\ ... \\ v^{n-1} \\ v^n
\end{array}
\right) ,
\end{eqnarray}
where $\lambda=\lambda(s)$ is some function.

Besides, the vector-function $v(s)$ has to be of unit length:
\begin{eqnarray}\label{ode_length}
\sum\limits_{j=1}^n (v^j)^2\equiv 1.
\end{eqnarray}

From (\ref{ode_sys}) we get
\begin{eqnarray*}
\frac{d}{ds} \left( \sum\limits_{j=1}^n (v^j)^2 \right) = 2 \left( -v^1 +\lambda \sum\limits_{j=1}^n (v^j)^2  \right) .
\end{eqnarray*}
Together with (\ref{ode_length}) this implies that $\lambda=v^1$, otherwise (\ref{ode_sys})-(\ref{ode_length}) has no solutions. Moreover, in this case we have
\begin{eqnarray*}
\frac{d}{ds} \left( \sum\limits_{j=1}^n (v^j)^2 - 1 \right) = 2 v ^1 \left( \sum\limits_{j=1}^n (v^j)^2 -1 \right) ,
\end{eqnarray*}
hence it is easy to see that (\ref{ode_length}) follows from (\ref{ode_sys}) if the initial data for (\ref{ode_sys}) satisfy
\begin{eqnarray}\label{ode_length0}
\sum\limits_{j=1}^n (v^j(0))^2 - 1 = 0.
\end{eqnarray}

Thus, for reconstructing the tractrix $\gamma$ from the given directrix $\tilde\gamma$, we have to solve the following Cauchy problem:
\begin{eqnarray}\label{ode_sys1}
\frac{d}{ds}
\left(
\begin{array}{l}
v^1 \\ v^2 \\ ... \\ v^{n-1} \\ v^n
\end{array}
\right) =
\left(
\begin{array}{lllllll}
0&k_1&  ...  & 0 & 0 \\
-k_1 & 0 &  ...  & 0 & 0 \\
...& ... & ...  & ... & ... \\
0 & 0 &  ...  & 0 & k_{n-1} \\
0 & 0 &  ...  &  -k_{n-1} & 0 \\
\end{array}
\right) \cdot
\left(
\begin{array}{l}
v^1 \\ v^2 \\ ... \\ v^{n-1} \\ v^n\\
\end{array}
\right)
 +
 \left( \begin{array}{l}
(v^1)^2 - 1\\ v^1v^2 \\ ... \\ v^1v^{n-1} \\ v^1v^n
\end{array}
\right) ,
\end{eqnarray}
\begin{eqnarray}\label{ode_data1}
\left(
\begin{array}{l}
v^1 \\ v^2 \\ ... \\ v^{n-1} \\ v^n\\
\end{array}
\right)_{s=0} \, = \,
\left(
\begin{array}{l}
v^1_0 \\ v^2_0 \\ ... \\ v^{n-1}_0 \\ v^n_0\\
\end{array}
\right) , \quad\quad \sum\limits_{j=1}^n (v^j_0)^2\equiv 1,
\end{eqnarray}
where $k_1(s)$, ... , $k_{n-1}(s)$ are given by the directrix $\tilde \gamma$, whereas $v^1(s)$, ..., $v^n(s)$ determine the desired tractrix $\gamma$.

Evidently, the system of first order ode's (\ref{ode_sys1}) with arbitrary $k_1(s)$, ..., $k_{n-1}(s)$ is of independent interest for studying from the analytical point of view. In the next section we will consider the particular case of $n=3$ with $k_1\equiv const$, $k_2\equiv 0$, which describes circular tractrices in $E^3$ whose directrices are circles.

Notice that  one has
\begin{eqnarray*}
 \frac{d f}{ds} = v^1 \left( v^1 \xi_1 + ... + v^n \xi_n \right), \quad\quad \left\vert \frac{d f}{ds} \right\vert = \left\vert v^1 \right\vert
\end{eqnarray*}
because of (\ref{tan}) and (\ref{ode_sys1}). Hence, the reconstructed tractrix $\gamma$ is regular everywhere except points where $v^1$ vanishes. Geometrically, the vanishing of $v^1$ means that $v$ becomes orthogonal to $\tilde\gamma$, c.f. \cite{tabachnikov}.

\bigskip

\bigskip

\textbf{3. Circular tractrices in $E^3$}

\medskip

Suppose $\tilde \gamma$ is a circle of radius $r$ in $E^3$. In appropriate Cartesian coordinates in $E^3$ the curve $\tilde\gamma$ is represented by the position vector
\begin{eqnarray}\label{pvector}
\tilde f (s) \, = \,
\left(
\begin{array}{l}
r\cos\frac{s}{r} \\ r\sin\frac{s}{r} \\  0 \\
\end{array}
\right),
\end{eqnarray}
its Frenet frame is given by
\begin{eqnarray}\label{FFrame}
\xi_1 (s) \, = \,
\left(
\begin{array}{l}
-\sin\frac{s}{r} \\ \cos\frac{s}{r} \\  0 \\
\end{array}
\right),  \quad
\xi_2 (s) \, = \,
\left(
\begin{array}{l}
-\cos\frac{s}{r} \\ -\sin\frac{s}{r} \\  0 \\
\end{array}
\right),  \quad
\xi_3 (s) \, = \,
\left(
\begin{array}{l}
0 \\ 0 \\  1 \\
\end{array}
\right),
\end{eqnarray}
and the curvatures are
\begin{eqnarray}\label{curv}
k_1 \equiv \frac{1}{r},\quad \quad k_2\equiv 0.
\end{eqnarray}

For reconstructing a circular tractrix $\gamma$ in $E^3$, whose directrix is the given circle $\tilde\gamma$, we have to solve the following system of ode's which is a particular case of (\ref{ode_sys1}),
\begin{eqnarray}\label{ode_sys2}
\left\{
\begin{array}{l}
\displaystyle{\frac{dv^1}{ds}} = \displaystyle{\frac{1}{r}} v^2 + (v^1)^2 - 1 \\ \\ \displaystyle{\frac{dv^2}{ds}} = - \displaystyle{\frac{1}{r}} v^1 + v^1v^2 \\ \\ \displaystyle{\frac{dv^3}{ds}} = v^1v^3
\end{array}
\right.
\end{eqnarray}
together with the equation
\begin{eqnarray}\label{ode_sys3}
(v^1)^2+(v^2)^2+(v^3)^2 = 1.
\end{eqnarray}

The solvability of (\ref{ode_sys2})-(\ref{ode_sys3}) strongly depends on whether $r$ is greater, equal or less than 1.

\bigskip

{\bf Theorem 1.}

\medskip

{\it 1. If $r>1$ then the solution of (\ref{ode_sys2}) is either
\begin{eqnarray}\label{large1}
v^1 = \displaystyle{-\frac{\lambda\, \sinh (\lambda s + c_3)}{\frac{c_1}{r}+\cosh (\lambda s + c_3)}}, \quad v^2 = \displaystyle{\frac{c_1+\frac{1}{r}\, \cosh (\lambda s + c_3)}{\frac{c_1}{r}+\cosh (\lambda s + c_3)}}, \quad v^3 = \displaystyle{\frac{\lambda\, c_2}{\frac{c_1}{r}+\cosh (\lambda s + c_3) }},
\end{eqnarray}
where $\lambda=\frac{\sqrt{r^2-1}}{r}$ and $c_1$, $c_2$, $c_3$ are arbitrary constants subject to $c_1^2+c_2^2 =1$,
or
\begin{eqnarray}\label{large2}
v^1 = \pm \displaystyle{\frac{\sqrt{r^2-1}}{r}}, \quad v^2 = \displaystyle{\frac{1}{r}}, \quad v^3 = 0.
\end{eqnarray}

2. If $r=1$ then the solution of (\ref{ode_sys2}) is either
\begin{eqnarray}\label{unit1}
v^1 = \displaystyle{-\frac{2s+c_3}{c_1+(s+c_3)^2}}, \quad v^2 = \displaystyle{1 - \frac{2}{c_1+(s+c_3)^2}}, \quad v^3 = \displaystyle{\frac{2 c_2}{c_1+(s+c_3)^2}},
\end{eqnarray}
where $c_1$, $c_2$, $c_3$ are arbitrary constants subject to $c_1 = 1+ c_2^2$,
or
\begin{eqnarray}\label{unit2}
v^1 = 0, \quad v^2 = 1, \quad v^3 = 0.
\end{eqnarray}

3. If $r<1$ then the solution of (\ref{ode_sys2}) is either
\begin{eqnarray}\label{small1}
v^1 = \displaystyle{\frac{\lambda\, \sin(\lambda s + c_3)}{\frac{c_1}{r}+\cos(\lambda s + c_3)}}, \quad v^2 = \displaystyle{\frac{c_1+\frac{1}{r} \cos(\lambda s + c_3)}{\frac{c_1}{r}+\cos(\lambda s + c_3)}}, \quad v^3 = \displaystyle{\frac{\lambda\, c_2}{\frac{c_1}{r}+\cos(\lambda s + c_3)}},
\end{eqnarray}
where $\lambda=\frac{\sqrt{1-r^2}}{r}$ and $c_1$, $c_2$, $c_3$ are arbitrary constants subject to $c_1^2-c_2^2 =1$,
or
\begin{eqnarray}\label{small2}
v^1 = 0, \quad v^2 = r, \quad v^3 = \pm \sqrt{1-r^2}.
\end{eqnarray} }

The proof is based on a detailed analysis of the system in question by elementary methods of the theory of ordinary differential equations. Namely, $v^2(s)$ is explicitly defined from the first equation of (\ref{ode_sys2}), and then the second equations rewrites as a second order ode for the single function $v^1(s)$:
\begin{eqnarray*}
\frac{d^2v^1}{ds^2}-3v^1\frac{dv^1}{ds} + v^1\left( \frac{1}{r^2}-1 \right) + (v^1)^3 = 0.
\end{eqnarray*}
Finally, $v^3(s)$ is explicitly determined up to a constant factor from the third equation of (\ref{ode_sys2}).

Consequently, if one substitutes the found solutions into (\ref{1}) keeping in mind (\ref{3}) and (\ref{pvector})-(\ref{FFrame}), then one obtains the position vector $f(s)$ of the desired circular tractrix $\gamma$ whose directrix is the given circle $\tilde\gamma$:
\begin{eqnarray}\label{finalTrac}
f(s) = \left(
\begin{array}{c}
r\cos{\frac{s}{r}} - v^1 \sin{\frac{s}{r}} - v^2 \cos{\frac{s}{r}} \\
r\sin{\frac{s}{r}} + v^1 \cos{\frac{s}{r}} - v^2 \sin{\frac{s}{r}} \\
v^3
\end{array}
\right) , \quad \quad s\in \mathbb R.
\end{eqnarray}

In the case of general solutions (\ref{large1}), (\ref{unit1}), (\ref{small1}), the vector-function $f(s)$ depends on three additional parameters $c_1$, $c_2$ and $c_3$ subject to one constrain depending on $r$. Clearly, the choice of $c_3$ corresponds to the rotations around the $x^3$-axe  (along the circle $\tilde\gamma$) in $E^3$, whereas $c_1$ and $c_2$ exert essential influence on the extrinsic form of the circular tractrix $\gamma$. Thus, up to rotations in $E^3$, we have a one-parametric family of non-congruent circular tractrices with the same directrix $\tilde\gamma$.

The exceptional solutions (\ref{large2}), (\ref{unit2}), (\ref{small2}), which are just the constant (stationary) solutions of (\ref{ode_sys2})-(\ref{ode_sys3}), correspond to exceptional circular tractrices. In the case of (\ref{large2}) one gets that $\gamma$ is the circle of radius $\sqrt{r^2-1}$ concentric to $\tilde \gamma$ and situated in the same horizontal coordinate plane as $\tilde \gamma$. In the case of either (\ref{unit2}) or (\ref{small2}) the corresponding circular tractrix $\gamma$ degenerates to a point.

\bigskip

\textbf{4. Circular Dini surfaces in $E^4$}

\medskip

Now view the circular tractrix $\gamma$ represented by (\ref{finalTrac}) as a curve in $E^4$ by adding to its position vector the fourth component equal to zero, and consider a two-dimensional surface $F^2$ in $E^4$ represented by the position vector
\begin{eqnarray}\label{GDini}
f(s,\varphi) =
\left(
\begin{array}{cccc}
\cos a\varphi& - \sin a\varphi& 0& 0 \\
\sin a\varphi& \cos a\varphi& 0& 0 \\
0&0&\cos b\varphi& - \sin b\varphi \\
0&0& \sin b\varphi& \cos b\varphi
\end{array}
\right)
\cdot
\left(
\begin{array}{c}
r\cos{\frac{s}{r}} - v^1 \sin{\frac{s}{r}} - v^2 \cos{\frac{s}{r}} \\
r\sin{\frac{s}{r}} + v^1 \cos{\frac{s}{r}} - v^2 \sin{\frac{s}{r}} \\
v^3 \\
0
\end{array}
\right) ,
\end{eqnarray}
where $a$, $b$ are arbitrary non-zero constants, and $(s,\varphi) \in \mathbb R^2$.

Geometrically, the surface $F^2$ is obtained by applying to the circular tractrix $\gamma$ a skew rotations in $E^4$ composed of two rotations in mutually orthogonal two-dimensional subspaces of $E^4$. Thus, $F^2$ is foliated by mutually congruent copies of $\gamma$ which we view as meridians of $F^2$.

\bigskip

{\bf Theorem 2.}

\medskip

{\it

In the case of $r>1$, if $\gamma$ results from (\ref{large1})
then the Gauss curvature $K$ of  $F^2$ is equal to
\begin{eqnarray*}
K\equiv \frac{c_2^2 (a^2-b^2)}{a^2(r^2-1)+b^2c_2^2}.
\end{eqnarray*}

In the case of $r=1$, if  $\gamma$ results from (\ref{unit1})
then the Gauss curvature $K$ of  $F^2$ is equal to
\begin{eqnarray*}
K\equiv \frac{c_2^2 (a^2-b^2)}{a^2+b^2c_2^2}.
\end{eqnarray*}

In the case of $r<1$, if  $\gamma$ results from (\ref{small1})
then the Gauss curvature $K$ of  $F^2$ is equal to
\begin{eqnarray*}
K\equiv \frac{c_2^2 (a^2-b^2)}{a^2(1-r^2)+b^2c_2^2}.
\end{eqnarray*}

Thus, in all three cases the Gauss curvature is constant.
}

\medskip

The proof is based on the use of standard formulae from the  elementary differential geometry concerning the first fundamental form, Christoffel symbols and Gauss curvature. We just notice that the first fundamental form $g$ of $F^2$ reads as follows:

if $r>1$, then
\begin{eqnarray*}
g=\displaystyle{\frac{\lambda^2\, \sinh^2(\lambda s + c_3)}{(\frac{c_1}{r}+\cosh (\lambda s + c_3))^2}}\cdot (ds + r\alpha d\varphi)^2 + \displaystyle{\frac{\lambda^2(\alpha^2(r^2-1)+\beta^2c_2^2)}{(\frac{c_1}{r}+\cosh (\lambda s + c_3))^2}}\cdot (d\varphi)^2;
\end{eqnarray*}

if $r=1$, then
\begin{eqnarray*}
g=\displaystyle{\frac{4\,(s + c_3)^2}{(c_1+(s + c_3)^2)^2}}\cdot (ds + \alpha d\varphi)^2 + \displaystyle{\frac{4(\alpha^2+\beta^2c_2^2)}{(c_1+(s + c_3)^2)^2}}\cdot (d\varphi)^2;
\end{eqnarray*}

if $r<1$, then
\begin{eqnarray*}
g=\displaystyle{\frac{\lambda^2\, \sin^2(\lambda s + c_3)}{(\frac{c_1}{r}+\cos (\lambda s + c_3))^2}}\cdot (ds + r\alpha d\varphi)^2 + \displaystyle{\frac{\lambda^2(\alpha^2(1-r^2)+\beta^2c_2^2)}{(\frac{c_1}{r}+\cos (\lambda s + c_3))^2}}\cdot (d\varphi)^2.
\end{eqnarray*}

\bigskip

\bigskip

\textbf{5. Open problems}

\medskip

1. To explore intrinsic and extrinsic geometrical properties of the circular Dini surfaces in $E^4$. For instance, it would be interesting to identify domains of the Euclidean / spherical / hyperbolic spaces which are isometrically realized on regular parts of the circular Dini surfaces in $E^4$.

2. To explicitly describe or classify {\it helicoidal} tractrices in $E^n$, $n\geq 3$, whose directrices are curves of constant curvatures, and then verify whether that helicoidal tractrices can be used for constructing novel examples of rotationally invariant surfaces / submanifolds of constant curvature in $E^{n+m}$, $n+m\geq 4$, which inherit geometric properties of the Beltrami and Dini surfaces related to the degeneracy of their Bianchi-Backlund transformations.

3. To examine the solvability of (\ref{ode_sys1}) with particular attention to qualitative and asymptotic properties of solutions. Ever the case of constant $k_1$, ..., $k_{n-1}$   is of interest.

\end{document}